\documentclass[]{amsart}
\usepackage{amsmath,amsfonts,amssymb}
\usepackage[english]{babel}
\def\N{\mathbb N}
\def\R{\mathbb R}

\def\C{\mathbb C}

\def\diff{{\rm Diff}}

\def\Cd#1{(\C^#1,0)}

\def\coef{\mathop{\rm Coef}\nolimits}
\def\Exp{\mathop{\rm Exp}\nolimits}

\theoremstyle{plain}
\newtheorem{theorem}{Theorem}[section]
\newtheorem{lemma}[theorem]{Lemma}

\newtheorem{proposition}[theorem]{Proposition}
\def\proof{{\it Proof: }}
\def\qed{\hfill\hbox{$\square$}}

\theoremstyle{definition}

\author[F.E. Brochero Mart\'{\i}nez]{F. E. Brochero Mart\'{\i}nez}
\address{
Departamento de Matem\'atica\\
UFMG\\
Belo Horizonte, MG\\
 30123-970\\
 Brazil\\
 }
 \email{fbrocher@mat.ufmg.br }
\thanks{The first author was  supported by CAPES, Brazil, Process:
BEX3083/05-5}

\author[L. L\'opez-Hernanz]{L. L\'opez-Hernanz}
\address{Departamento de \'Algebra, Geometr\'{\i}a y Topolog\'{\i}a\\
Universidad de Valladolid, Spain} \email{llopez@agt.uva.es}
\thanks{The second author was  supported by FPU, Spain, Process: AP2005-3784}

\date{\today
}


\title{Class of the infinitesimal generator  of a diffeomorphism in $\Cd m$}

\begin{document}

\begin{abstract}
Let $F$ be an analytic diffeomorphism in $\Cd m$ tangent to the
identity of order $n$. The infinitesimal generator of $F$ is the
formal vector field $X$ such that $\Exp X=F$. In this paper we
provide an elementary proof of the fact that $X$ belongs to the
Gevrey class of order $1/n$.
\end{abstract}

\maketitle

\section{Introduction}
For each couple of integers $m\ge 1$ and $n\ge 2$, let us denote
${\hat{\mathfrak X}}_{n}\Cd m$ the module of formal vector fields
of order  $\ge n$ in $\Cd m$  and $\widehat\diff_{n}\Cd m$ the
group of formal diffeomorphisms in $\Cd m$ tangent to the identity
of order $\ge n$, i.e, $F\in \widehat\diff_{n}\Cd m$ if and only
if $ \nu(F):=\min\{\nu_0(x_i\circ F-x_i)|i=1,\dots,m\}-1\ge n$.
 For any  $X\in {\widehat{\mathfrak X}}_{n}\Cd m$, the
exponential operator of $X$ is  the application $\exp{X}:\C[[x]]\to
\C[[x]]$ defined by the formula
$$\exp{X}(g)=\sum_{j=0}^\infty \frac
1{j!}X^j(g)$$ where $X^0(g)=g$ and $X^{j+1}(g)=X(X^j(g))$.  It is a
classical result (for instance, see \cite{MaRa}) that the
application
$$\begin{array}{rcl}\Exp:{\widehat{\mathfrak X}}_{n}\Cd m&\to&
\widehat\diff_{n-1}\Cd m\\
X&\mapsto&(\exp X(x_1),\dots,\exp X(x_m))
\end{array}$$
is a bijection. The formal vector field $X$ such that $F=\Exp(X)$ is
called the {\em infinitesimal generator} of $F$.

Let $x=(x_1,\dots, x_m)$ and for any $s\in \R$ let $\C[[x]]_s$
denote the subset of elements of $\C[[x]]$ that satisfy the
$s$-Gevrey condition, i.e.
$$f(x)=\sum_{k=0}^\infty f_k(x)\in \C[[x]]_s\quad\hbox{if and only
if}\quad  \sum_{k=0}^\infty \frac{f_k(x)}{k!^s}\in \C\{x\},$$
where $f_k(x)$ is  homogeneous of degree $k$. Let us observe that
0-Gevrey condition means analyticity, and
$\C\{x\}\subset\C[[x]]_s\subset \C[[x]]_t$ if $0<s<t$. Let
$\mathfrak X_{n}\Cd m_s\subseteq \hat{\mathfrak X}_{n}\Cd m$ be
the set of $s$-Gevrey vector fields $X=\sum\limits_{k=1}^m
X(x_k)\frac {\partial}{\partial x_k}$ with $X(x_k)\in \C[[x]]_s$
and $\diff_{n}\Cd m_s= \widehat{\diff}_{n}\Cd m\cap (\C[[x]]_s)^m$
the set of  $s$-Gevrey diffeomorphisms  tangent to the identity of
order $\ge n$.

We will prove the following result
\begin{theorem}\label{biyeccion} For any $s\ge \frac 1 {n-1}$ the
application $\Exp$ gives a bijection
$$\Exp:\mathfrak X_{n}\Cd m_s\to \diff_{n-1}\Cd m_s.$$
 In particular, the infinitesimal generator of any tangent to the
identity  analytic diffeomorphism $F$ is $\frac{1}{\nu (F)}$-Gevrey.
\end{theorem}

 In general,
 $X$ may be  divergent for a convergent $F$, for instance, Szekeres
\cite{Sze} and Baker \cite{Bak} proved that every entire
holomorphic function tangent to the identity of order $k$ in
dimension $1$ has a non-convergent infinitesimal generator,  Ahern
and Rosay \cite{AhRo} proved that this kind of diffeomorphisms
cannot be the time-1 map of a $C^{3k+3}$-vector field, and finally
J. Rey \cite{Rey} showed that it cannot be the time-1 map of a
$C^{k+1}$-vector field, which is the best possible bound. Thus,
the map $\Exp:\mathfrak X_n\Cd m_0\to \diff_{n-1}\Cd m_0$ is not
surjective for any couple of positive integers $m,n$. In addition,
in dimension $1$, using resummation arguments, it is proved that
if an analytic diffeomorphism $f(x)=x+a_{k+1}x^{k+1}+\cdots$ with
$a_{k+1}\ne 0$ has a divergent infinitesimal generator $X$, then
$X$ is $k$-summable, so $X$ is Gevrey of order $\frac 1k$, but not
smaller
 (see \cite{MaRa2}, \cite{Eca} and \cite{Mal}).
Therefore, the condition $s\ge \frac 1{n-1}$ is necessary.
\vspace{2mm}

{\em Acknowledgements.} We would like to
 thank Javier Rib\'on for pointing out the above results on convergence and resummation
 in dimension 1, and giving us
 the idea to improve our first version of this paper.
We would like to thank Jos\'e Cano for fruitful conversation and
computational calculations and Felipe Cano for useful
conversation.

\section{Technical estimations}
In this paper, we take the following notations:
\begin{itemize}
 \item $h_k(x)$ will denote
the homogeneous polynomial $\sum\limits_{\alpha\in \N^m \atop
|\alpha |=k} x^\alpha$.\vspace{-3mm}
\item $H_{s,n}(x)$ the series $\sum\limits_{q=n}^\infty (q+m-n)!^s
h_q(x)$.
\item ${\frac{\partial}{\partial x}}$ the differential operator
$\sum\limits_{k=1}^m \frac {\partial}{\partial x_k}$.
\end{itemize}
For formal series  $f(x)=\sum_{\alpha} f_\alpha x^\alpha$ and
$g(x)=\sum_{\alpha} g_\alpha x^\alpha$,  we say that $f\preceq g$ if
$|f_\alpha|\le |g_\alpha|$ for any $\alpha\in \N^m$. We get in this
way a partial order in $\C[[x]]$, and also in
 ${\hat{\mathfrak X}}_n\Cd m$
and $\widehat{\diff}_n\Cd m$, working on the component function.
From the definition of Gevrey condition, it can be seen that $X\in
\mathfrak X_n\Cd m_s$ if and only if there exists $a\in \R^+$ such
that
$$\coef _q(X)\preceq (q+m-n)!^s a^q h_q(x){\frac{\partial}{\partial x}},$$
where $\coef _q(X)$ denotes the homogeneous term of $X$ of degree
$q$. Thus $X\in\mathfrak X_n\Cd m_s$ if and only if there exists
$a\in \R^+$ such that $X\preceq H_{s,n}(ax){\frac{\partial}{\partial
x}}.$

We need the following technical lemmas:

\begin{lemma}For every $k,l\in \N^*$
$$h_k{\frac{\partial}{\partial x}} h_l\preceq
(l+m-1)\min\Bigl\{\binom{k+m-1}{m-1},\binom{l+m-2}{m-1}\Bigl\}
h_{k+l-1}.$$
\end{lemma}
\proof Observe that
\begin{eqnarray*}
{\frac{\partial}{\partial x}} h_l&=&\sum_{k=1}^m \frac
{\partial}{\partial x_k} \sum_{\alpha\in \N^m \atop |\alpha|=l}
x^\alpha =\sum\limits_{k=1}^m
 \sum\limits_{\alpha\in \N^m \atop
|\alpha|=l} \alpha_k
\frac {x^{\alpha}}{x_k}\\
&=&\sum\limits_{\beta\in \N^m \atop |\beta|=l-1}\sum\limits_{k=1}^m
(\beta_k+1) x^\beta=(l+m-1)h_{l-1}
\end{eqnarray*}
Now,  the coefficient of $x^\alpha$ in the product
$h_k(x)h_{l-1}(x)$
is less than or equal to  the minimum  between the number of
monomials  of $h_k$ and the number of monomials of $h_{l-1}$, and
the number of monomials of $h_j$ is $\binom{j+m-1}{m-1}$, that
corresponds to the number of ordered partitions of $j$ in $m$
parts; therefore,
$$h_k{\frac{\partial}{\partial x}} h_l=(l+m-1)h_kh_{l-1}\preceq
(l+m-1)\binom{\min\{k,l-1\}+m-1}{m-1} h_{k+l-1}. \eqno\qed$$

\begin{lemma}\label{theta}Let $\Theta(y)=\sum\limits_{j=n}^\infty \binom{m-1+j}{m-1}y^{j-n}$.
Then $\Theta(y)$ converges for any $|y|<1$.
\end{lemma}
\proof  Since $\sum\limits_{j=n}^\infty  y^{m-1+j}=\frac
{y^{m+n-1}}{1-y}$ converges for any $|y|<1$ then
$$\Theta(y)=\frac 1{(m-1)!}\frac 1{y^n}\frac {d^{m-1}}{dy^{m-1}}\left(\frac
{y^{m+n-1}}{1-y}\right)$$ converges for any $|y|<1$.\qed

\begin{lemma}\label{acotado} For any $s>0$ and integers $m\ge 1$ and  $n\ge 2$, the
sequence $\{b_q\}_{q\ge 2n-1}$ given by
$$b_q=\sum\limits_{j=n}^{\lfloor \frac
{q+1}2\rfloor}\left(\frac{(j+m-n)!(q-j+1+m-n)!}{m!(q+m-n)!}(q-j+m)^{n-1}\right)^s\binom{j+m-1}{m-1},$$
is bounded.
\end{lemma}

\proof Observe that
\begin{eqnarray*}
\frac {(q-j+m)^{n-1}}{(q-j+2+m-n)\cdots
(q-j+m)}&<&\left(\frac{q-j+m}{q-j+2+m-n}\right)^{n-1}\\
&\le&\left(\frac{\frac{q-1}{2}+m}{\frac{q-1}{2}+2+m-n}\right)^{n-1}\le\left(\frac{m+n-1}{m+1}\right)^{n-1}
\end{eqnarray*}
 then
$$b_n\le\left(\frac{m+n-1}{m+1}\right)^{s(n-1)}\sum\limits_{j=n}^{\lfloor \frac
{q+1}2\rfloor}\left(\frac{(j+m-n)!(q-j+m)!}{m!(q+m-n)!}\right)^s\binom{j+m-1}{m-1}.$$
In addition
$$\frac{m+1}{q+m-j+1}<\frac{m+2}{q+m-j+2}<\cdots<\frac{j+m-n}{q+m-n}$$
and
$$\frac{j+m-n}{q+m-n}\le
\frac{\frac{q+1}{2}+m-n}{q+m-n}\le \max \left\{\frac{1}{2},
\frac{m}{m+n-1}\right\}= C_{m,n}<1;$$ from lemma \ref{theta},
$$b_q<
\left(\frac{m+n-1}{m+1}\right)^{s(n-1)}\Theta ( C^s_{m,n}).\eqno
\qed$$

\begin{proposition}\label{potencias} Let $s\ge \frac 1
{n-1}$,
 $ X\in \widehat{\mathfrak X}_n\Cd m$ and $a\in \R^+$ such that $$\coef _q(X)\preceq
 (q+m-n)!^s
 a^q h_q(x)\frac{\partial}{\partial x}$$ for all $n\le q\le N$, and let us denote
 $A=2m!^s \left(\frac{m+n-1}{m+1}\right)^{s(n-1)}
 \Theta(C_{m,n}^s)$. For every
$q,k$ with $n\le q\le N+k-1$,
$$\coef _q(X^k)\preceq (a A)^{k-1}(q+m-n)!^s
 a^q h_q(x)\frac{\partial}{\partial x},$$
\end{proposition}

\proof
Since $X^k=\sum\limits_{i=1}^m X^k(x_i)\frac{\partial}{\partial
x_i}$, it is enough to prove the affirmation for $X^k(x_i)$, where
$i\in \{1,2,\dots,m\}$. Let us write $X=\sum\limits_{j=n}^\infty
X_j$, where $X_j$ is homogeneous of degree $j$. We will proceed by
induction on $k$; if $k=1$, by hypothesis
$$X_q(x_i)\preceq (q+m-n)!^sa^q h_q(x)\quad\hbox{for every $n\le q\le N$.}$$
Suppose that the lemma is true for every $k\le p$, then, since the
order of $X^j$ is greater than or equal to $(n-1)j+1$,
$\coef_q(X^{p+1})=0$ for $2\le q\le (n-1)p+n-1$ and for
$(n-1)p+n\le q\le N+p$ we have {\small
$$
\begin{array}{l}
\phantom{i}\hspace{-5mm} \coef
_q(X^{p+1}(x_i))=\coef_q(X(X^p(x_i)))=
 \coef _q\Bigl(\sum\limits_{j=n}^{\infty}X_j(X^p(x_i))\Bigr)\\
=\sum\limits_{j=n}^{q-(n-1)p} X_j \coef _{q+1-j}(X^p(x_i))\\
\preceq\sum\limits_{j=n}^{q-(n-1)p} (j+m-n)!^s a^j
h_j(x){\frac{\partial}{\partial x}}\Bigl((aA)^{p-1}(q-j+1+m-n)!^s  a^{q+1-j} h_{q+1-j}(x)\Bigr)\\
\preceq \sum\limits_{j=n}^{q-n+1}(j+m-n)!^s(q-j+1+m-n)!^s
(q-j+m)\binom{\min\{j,q-j\}+m-1}{m-1}
A^{p-1}a^{q+p} h_q,\\
\preceq 2 \sum\limits_{j=n}^{\lfloor \frac
{q+1}2\rfloor}((j+m-n)!(q-j+1+m-n)!(q+m-j)^{n-1})^s\binom{j+m-1}{m-1}
A^{p-1}a^{q+p} h_q.
\end{array}$$}%
Now, observe that  $$ b_jm!^s(q+m-n)!^s=\sum\limits_{j=n}^{\lfloor
\frac
{q+1}2\rfloor}((j+m-n)!(q-j+1+m-n)!(q-j+m)^{n-1})^s\binom{j+m-1}{m-1},$$
where $\{b_q\}$  is the sequence defined in lemma \ref{acotado}; it
follows that
\begin{eqnarray*}\coef_q(X^{p+1}(x_i))&\preceq& 2  b_q m!^s(q+m-n)!^s
A^{p-1}a^{q+p} h_q\\
&\preceq& (q+m-n)!^s (aA)^p a^q h_q\qquad \qed
\end{eqnarray*}

\section{Proof of theorem \ref{biyeccion}.}
To prove that the application $\Exp:\mathfrak X_{n}\Cd m_s\to
\diff_{n-1}\Cd m_s$ is well defined for $s\ge \frac 1 {n-1}$, let
$X\in \mathfrak{X}_n \Cd m_s$, $a>0$ be such that $X\preceq
H_{s,n}(ax)$, and $A$  as in proposition \ref{potencias}.
Then by proposition \ref{potencias} we have
\begin{eqnarray*}
\coef _q(\exp X(x_j))&=&\sum_{k=1}^\infty \frac 1{k!}\coef _q(X^k(x_j))\\
&\preceq&  \sum_{k=1}^\infty \frac 1{k!} (aA)^{k-1} (q+m-n)!^s
a^q h_q(x)\\
\end{eqnarray*}
therefore $\Exp(X)\preceq  \sum\limits_{k=1}^\infty \frac
{(aA)}{k!}\!{\,\!^{k-1}\atop \ }
 H_{s,n}(ax)$. Now, to prove that $\Exp$ is
 surjective, let us consider a diffeomorphism
$F(x)=(x_1+f_1(x),\dots ,x_m+f_m(x))\in \diff_{n-1}\Cd m_s$ where
$f_j(x)=\sum\limits_{q=n}^\infty f_{j,q}(x)\in \C[[x]]_s$ and
$f_{j,q}(x)$ is an homogeneous polynomial of degree $q$. Then
there exists $a>0$ such that $f_{j, q}(x) \preceq (q+m-n)!^s a^q
h_q(x)$. Observe that, making a linear change of coordinates, we
can suppose that $a$ is small enough such that
$\sum\limits_{k=2}^\infty \frac 1{k!} (2aA)^{k-1}\le\frac 12$. If
$X=\sum\limits_{q=n}^\infty X_q$ is the infinitesimal generator of
$F(x)$, we will by induction on $q$ that
$$X_q\preceq (q+m-n)!^s (2a)^q h_q(x){\frac{\partial}{\partial x}}
.$$ For $q=n$
$$X_n(x_j)=f_{j,n}(x)\preceq m!^s a^n h_n(x)\preceq m!^s (2a)^n
h_n(x).$$ Suppose that the claim is true for  any integer between
$n$ and $q$, it follows that
$$f_{j, q+1}(x)=\coef _{q+1}\Bigl(\sum_{k=1}^\infty \frac 1{k!}
X^k(x_j)\Bigr)=X_{q+1}(x_j)+\sum_{k=2}^q \frac 1{k!} {\mathcal
Coef}_{q+1}\bigl(X^k(x_j)\bigr),$$
 using  proposition
\ref{potencias}
\begin{eqnarray*}
X_{q+1}(x_j)&\preceq& (q+1+m-n)!^s a^{q+1} h_{q+1}(x)\\
&&\hspace{2cm}+\sum_{k=2}^\infty \frac 1{k!} (2aA)^{k-1}
(q+1+m-n)!^s (2a)^{q+1} h_{q+1}(x)\\
&\preceq&\left( \frac 1{2^{q+1}}+\sum_{k=2}^\infty \frac 1{k!}
(2aA)^{k-1}\right)(q+1+m-n)!^s
(2a)^{q+1} h_{q+1}(x) \\
&\preceq&(q+1+m-n)!^s (2a)^{q+1} h_{q+1}(x),
\end{eqnarray*}
in other words $X\preceq H_{s,n}(2a)\frac{\partial}{\partial
x}$.\qed

\section{Case $0<s<\frac 1{n-1}$}
As we indicated in the introduction, in this case, there exists
$F\in \diff_{n-1} \Cd m_s$ such that its infinitesimal generator is
not $s$-Gevrey , but the reciprocal is true, i.e.

\begin{proposition}\label{biendefinido}Let $0\le s\le \frac1{n-1}$, and $X\in \mathfrak X_n \Cd m_{s}$. Then
$\Exp(X)\in \diff_{n-1}\Cd m_{s}$.
\end{proposition}

Observe that the case $s=0$ is a classical result about the
existence of solution of an analytic differential equation.  To
prove this proposition in the case $s>0$ we need the following
lemma

\begin{lemma}\label{radios} Let $t,r\in \R$ such that $0<t<1$ and $1-t<r<1$. Let  $\{a_k\}$ be the
sequence defined by $a_1=a>0$ and for $k\ge 1$,
$a_{k+1}=\sup\limits_{q\in \N^*}\sqrt[q+k]{\frac
{(q+m)^{1-t}}{(k+1)^r}} a_k$. Then $\{a_k\}$ is increasing and
convergent.
\end{lemma}

\proof Taking $q\gg k$ it clear that $\sqrt[q+k]{\frac
{(q+m)^{1-t}}{(k+1)^r}}>1$, and then $a_{k+1}>a_k$. Now, we know
by Bernoulli inequality that
$$
\sqrt[q+k]{\frac {q+m}{(k+1)^{\frac r{1-t}}}}<1+\frac
1{q+k}\left(\frac {q+m}{(k+1)^{\frac r{1- t}}}-1\right)
< 1+\frac 1{(k+1)^{\frac r{1-t}}}$$ for $k>m$, so
$$a_{k+1}<\left(1+\frac 1{(k+1)^{\frac r{1-t}}}\right)^{1-t}a_k<\Biggl(\prod_{j=m+1}^{k+1}
\Bigl(1+\frac 1{j^{\frac r{1-t}}}\Bigr)\Biggr)^{1-t}a_m,$$ and since
$\frac {r}{1-t}>1$ it follows that $\{a_k\}$ is bounded, thereby  it
is convergent.\qed

{\it Proof of proposition \ref{biendefinido}:}  If $s\in (0,\frac 1
{n-1})$,
 $ X\in {\mathfrak X}_n\Cd m_s$ and $a\in \R^+$ such that
 $X\preceq H_{s,n}(ax)\frac{\partial}{\partial x}$
 then  for $t=s(n-1)$, $r\in(1-t,1)$
   and $\{a_k\}$ as in lemma \ref{radios},
using the arguments of proposition \ref{potencias} and the fact
that
   $k^ra_k^{k+q-1}\ge (q+m)^{1-t} a_{k-1}^{k+q-1}$ for every $q\ge
   2$,  we can prove that
$$X^k\preceq (a_k A)^{k-1}k!^rH_{s,n}(a_kx)\frac{\partial}{\partial x},$$
 where
 $A=2m!^s \left(\frac{m+n-1}{m+1}\right)^{s(n-1)}
 \Theta(C_{m,n}^s)$.
Let  $c=\lim\limits_{k\to \infty} a_k$.
Therefore we have
$$
\coef_q(\exp(X)(x_j))=\sum_{k=1}^\infty \frac 1{k!}\coef_q(X^k(x_j))
\preceq \sum_{k=1}^\infty \frac {(cA)^{k-1}}{k!^{1-r}} (m+q-n)!^s
c^q h_q(x)
$$
Thus $\Exp(X)\preceq  \sum\limits_{k=1}^\infty \frac
{(cA)^{k-1}}{k!^{1-r}}
 H_{s,n}(cx)\frac{\partial}{\partial x}$.\qed

\end{document}